\documentclass[12pt]{amsart} \usepackage{amssymb} \unitlength=1pt
\newfont{\ams}{msbm10 at 12pt} \newfont{\amsi}{msbm8} 
\newcommand{\nc}{\newcommand} %\renewcommand{\u}{{\hbox{\bf u}}}
\nc{\C}{{\mathcal C}} \nc{\CB}{{\mathcal B}}
\nc{\supp}{\operatorname{supp}} \nc{\M}{{\mathcal M}}
\nc{\LL}{{\mathbf L}} \renewcommand{\P}{{\mathbf P}}
\renewcommand{\L}{{\mathcal L}} \nc{\D}{{\mathbb D}}
\renewcommand{\O}{{\mathcal O}} \nc{\St}{\operatorname{St}^{\bullet}}
\nc{\A}{{\frak A}} %\nc{\D}{{\frak D}} \renewcommand{\H}{{\frak H}}
\nc{\1}{{\mathbf 1}} \nc{\CC}{{\mathbb C}} \nc{\R}{{\mathbb R}}
\renewcommand{\k}{{\mathbb k}} \nc{\Q}{{\mathbb Q}} \nc{\U}{{\mathbf
U}} \nc{\B}{{\mathbf B}} \nc{\N}{{\mathbb N}} \nc{\Z}{{\mathbb Z}}
\nc{\F}{{\overline{\mathbb F}}} \nc{\CA}{{\mathcal A}}
\nc{\Rhom}{\operatorname{RHom}\bul} \nc{\Ad}{\operatorname{Ad}}
\nc{\Res}{\operatorname{Res}} \nc{\gr}{\operatorname{gr}}
\nc{\tr}{\operatorname{tr}} \nc{\End}{\operatorname{End}}

 \nc{\g}{{\frak g}}
\nc{\hatg}{\hat{\frak g}} 
\renewcommand{\b}{{\frak b}} \nc{\sem}{{$\S_{\g}^{\g_{>0}}}}
\nc{\gl}{{\frak g\frak l}} \nc{\n}{{\frak n}} \nc{\si}{{\frac\infty
2}} \nc{\p}{{\frak p}} \nc{\h}{{\frak h}} \renewcommand{\u}{{\frak u}}
\nc{\Ind}{\operatorname{Ind}} \nc{\ch}{\operatorname{ch}}
\nc{\Coind}{\operatorname{Coind}} \nc{\opp}{{\operatorname{opp}}}
\nc{\Ker}{\operatorname{Ker}} \nc{\im}{\operatorname{Im}}
\nc{\Coker}{\operatorname{Coker}}
\nc{\Cone}{\operatorname{Cone}}
\nc{\dirlim}{\underset{\rightarrow}{\operatorname{lim}}}
\nc{\invlim}{\underset{\leftarrow}{\operatorname{lim}}}
\nc{\Sem}{{{\mathbf S}_{\g}^{\g_{>0}}}} \nc{\CN}{{\mathcal N}}
\nc{\Ext}{\operatorname{Ext}^{\bullet}} \nc{\ext}{\operatorname{Ext}}
\nc{\tilW}{\til{W}}

\nc{\lth}{\ell t} \nc{\BB}{\mathcal{B}}
\nc{\Tor}{\operatorname{Tor}_{\bullet}} \nc{\tor}{\operatorname{Tor}}
\nc{\Tors}{\operatorname{Tor}_{\frac \infty 2+\bullet}}
\nc{\Exts}{\operatorname{Ext}^{\frac \infty 2+\bullet}}
\nc{\Hom}{\operatorname{Hom}^{\bullet}} \nc{\ad}{\operatorname{ad}}

\nc{\Spec}{\operatorname{Spec}}
\renewcommand{\hom}{\operatorname{Hom}}

\renewcommand{\mod}{\operatorname{-mod}} \nc{\Mod}{\operatorname{Mod}}

\nc{\Barb}{\operatorname{Bar}^{\bullet}}

\nc{\upX}{X^{\uparrow}} \nc{\upcD}{{\mathcal D}^{\uparrow}}
\nc{\upD}{D^{\uparrow}} \nc{\dX}{X^{\downarrow}} \nc{\dcD}{{\mathcal
D}^{\downarrow}} \nc{\dD}{D^{\downarrow}} \nc{\upC}{{\mathcal
C}^{\uparrow}} \nc{\dC}{{\mathcal C}^{\downarrow}}
\nc{\underA}{\underline{A}} \nc{\underC}{\underline{\CC}}
\nc{\underB}{\underline{B}} \nc{\underk}{\underline{\k}}
\nc{\Db}{D^{\bullet}} \nc{\ten}{{\otimes}}
\nc{\tenl}{\overset{\operatorname{L}}\ten} \nc{\map}{\longrightarrow}

\nc{\bs}{\bigskip\\} \nc{\ms}{\smallskip\\}

\nc{\tilbar}{\widetilde{\operatorname{Bar}}}
\nc{\tilBarb}{\widetilde{\operatorname{Bar}}^{\bullet}}
\nc{\overr}{\overline{R}} \nc{\overI}{\overline{I}}
\nc{\overX}{\overline{X}} \nc{\overh}{\overline{h}}
\nc{\overY}{\overline{Y}} \nc{\overW}{\overline{W}}
\nc{\linbar}{\overline{\operatorname{Bar}}}

\nc{\til}{\widetilde} \nc{\oppA}{A^{\sharp}} \nc{\Lemma}{{\bf Lemma:\
}} \nc{\Theorem}{{\bf Theorem:}\ } \nc{\Cor}{{\bf Corollary:}\ }
\nc{\Def}{{\bf Definition:}\ } \nc{\Prop}{{\bf Proposition:}\ }
\nc{\Con}{{\bf Conjecture:}\ } \nc{\Rem}{{\bf Remark:}\ }
\nc{\dok}{\noindent{\bf Proof.}\ }

\nc{\SInd}{\operatorname{S-Ind}} \nc{\SCoind}{\operatorname{S-Coind}}
\nc{\bul}{^{\bullet}} \nc{\stand}{C^{\frac\infty2+\bullet}}
\nc{\ssn}{\subsection{}} \nc{\sssn}{\subsubsection{}}

\nc{\hgt}{\operatorname{ht}} %\nc{\left}{\operatorname{left}}
\nc{\sqbinom}{\fracwithdelims[][0pt]}

\address{Independent University of Moscow, Pervomajskaya st.  16-18,
Moscow 105037, Russia} 
\email{hippie@@ium.ips.ras.ru} \author{Sergey
Arkhipov} 
\title{Semiinfinite cohomology of contragradient Weyl
modules over small quantum groups} 

\begin{document} 
\maketitle 
\section{Introduction}
The present paper can be considered as a natural extension the article
\cite{Ar7}.  Fix root data $(Y,X,\ldots)$ of the finite type
$(I,\cdot)$ and a positive integer number $\ell$.  In \cite{Ar7} we
obtained a nice description for semiinfinite cohomology of the trivial
module $\underC$ over the {\em small quantum group} $\u_\ell$
corresponding to the root data $(Y,X,\ldots)$ in terms of local
cohomology of the structure sheaf on the {\em nilpotent cone} $\CN$ in
the corresponding semisimple Lie algebra $\g$.

The geometric approach to various homological questions concerning the
algebra $\u_\ell$ appeared first in the pioneering paper of Ginzburg
and Kumar \cite{GK}.  Let us state the main result from that paper.

\vskip 1mm \noindent {\bf Theorem:}
$\Ext_{\u_\ell}(\underC,\underC)=H^0(\CN,\O_{\CN})$ as an associative
algebra.  The grading on the right hand side is provided by the action
of the group $\CC^*$ on the affine variety $\CN$.\qed

Moreover it was shown in \cite{GK} that both sides of the equality
carry natural structures of integrable $\g$-modules and the
isomorphism constructed in the paper takes the left hand side
$\g$-module structure to the right hand side one.

{\em Semiinfinite cohomology} of the trivial $\u_\ell$-module was
considered in \cite{Ar1}, \cite{Ar5}, \cite{Ar6} and \cite{Ar7}.
Consider semiinfinite cohomology of a graded associative algebra $A$
(see \cite{Ar1}, \cite{Ar2} for the definition of this cohomology
theory).  It is proved in \cite{Ar2} that the algebra
$\Ext_A(\underC,\underC)$ acts naturally on the semiinfinite
cohomology of $A$.  Thus in particular for a $\u_\ell$-module $M$ one
can treat $\Exts_{\u_\ell}(\underC,M)$ as a quasicoherent sheaf on
$\CN$.  The following statement known as the Feigin Conjecture was
proved in \cite{Ar7}.  Consider the standard positive nilpotent
subalgebra $\n^+\subset\CN\subset\g$.  \vskip 3mm \noindent \Theorem
The quasicoherent sheaf on $\CN$ provided by
$\Exts_{\u_\ell}(\underC,\underC)$ is equal to the sheaf of local
cohomology of $\O_{\CN}$ with support on $\n^+\subset\CN$.\qed

Moreover note that the simply connected Lie group $G$ with the Lie
algebra equal to $\g$ acts naturally on $\CN$.  This action provides a
structure of $\n^+$-integrable $\g$-modules on the local cohomology
spaces $H^{\sharp(R^+)}_{\n^+}(\CN,\O_{\CN})$.  On the other hand it
was shown in \cite{Ar1} that there exists a natural $U(\g)$-module
structure on $\Exts_{\u_\ell}(\underC,\underC)$.  It is proved in
\cite{Ar7} that the described $\g$-module structures coincide.

Consider the {\em contragradient Weyl module} $\D W_\ell(\ell\lambda)$
over $\u_\ell$ with the highest weight $\ell\lambda$.  In the present
paper we provide a geometric description of semiinfinite cohomology of
$\u_\ell$ with coefficients in $\D W_\ell(\ell\lambda)$.  To formulate
the exact statement we need some geometric notation.

Consider the flag variety $G/B$ for the group $G$ and its cotangent
bundle $T^*(G/B)$.  Below we denote these varieties by $\BB$ and
$\til{\CN}$ respectively.The natural projection $\til{\CN}\map\BB$ is
denoted by $p$.  The moment map for the symplectic $G$-action on
$\tilde{\CN}$ provides the {\em Springer-Grothendieck resolution} of
singularities of the nilpotent cone $\mu:\ \tilde{\CN}\map\CN$.
Consider the linear bundle $\L(\lambda)$ on $\BB$.  Note that by
\cite{GK} we have $$ \Ext_{\u_\ell}(\underC,\D
W_\ell(\ell\lambda))=H^0(\tilde{\CN},p^*\L(\lambda)).  $$ We call the
following statement the generalized Feigin Conjecture.  It is the main
result of the paper.

\vskip 3mm \noindent {\bf Conjecture:}  The
$H^0(\tilde{\CN},\O_{\tilde{\CN}})=\Ext_{\u_\ell}(\underC,\underC)$-module
$\Exts_{\u_\ell}(\underC,\D W_\ell(\ell\lambda))$ is isomorphic to
$H^{\sharp(R^+)}_{\mu^{-1}(\n^+)}(\tilde{\CN},p^*\L(\lambda))$.\qed

\ssn Let us describe briefly the structure of the paper.  In the
second section we recall in more detail results concerning both the
usual and the semiinfinite cohomology of small quantum groups
mentioned above.  Usung Ginsburg and Kumar's description of the
ordinary cohomology of $\u_\ell$ with coefficients in $\D
W_\ell(\ell\lambda)$ we construct a morphism of
$H^0(\tilde{\CN},\O_{\tilde{\CN}})$-modules $$ \sigma:\
H_{\mu^{-1}(\n^+)}^{\sharp(R^+)}(\tilde{\CN},p^*\L(\lambda))\map
\Exts_{\u_\ell}(\underC,\D W_\ell(\ell\lambda)).  $$

In the third section we construct a certain specialization of the
quantum BGG resolution into the root of unity.  We call the obtained
complex {\em the contragradient quasi-BGG complex} and denote it by
$\D B\bul_\ell(\mu)$.  This complex consists of direct sums of
{\em contragradient quasi-Verma modules} and its zero cohomology
module equals $\D W_\ell(\mu)$.  A quasi-Verma module $\D
M_\ell^w(w\cdot\mu)$ provides a certain specialization for the
family of the usual contragradient quantum Verma modules $\D
M_{\xi}(w\cdot)$ defined a priori at generic values of the
quantizing parameter $\xi$ into the root of unity.  However it turns
out that contragradient quasi-Verma modules {\em differ from the usual
contragradient Verma modules for the algebra} $\U_\ell$.  We show that
for a prime nubmber $\ell$ the contragradient quasi-BGG complex is
in fact quasiisomorphic to $\D W_\ell(\mu)$.  It is natural to
conjecture that the statement remains true for all but finitely many
roots of unity but our considerations do not provide the proof in the
general case.

 We perform the whole construction of the quasi-BGG complex over the
ring $\CA:=\Z[v,v^{-1},]$ and obtain a complex of $\U_{\CA}$-modules
$\D B\bul_{\CA}(\mu)$ for every integral dominant weight $\mu$.  Here
$\U_{\CA}$ denotes the Lusztig version of the quantum group with
quantum divided powers.  Next we recall the result Lusztig stating
that the specialization of the algebra $\U_{\CA}$ under the base
change from $\CA$ to $\F_\ell$ for $\ell$ simple coincides up to a
central extension with the specialization of the Kostant integral form
for the usual universal enveloping algebra $U_{\Z}(\g)$ to the same
field.  This way we obtain a geometric interpretation of the complex
$\D B\bul_{\F_\ell}(\mu)$ in terms of local cohomology of the flag
variety $\BB_{\F_\ell}$ in characteristic $\ell$ as follows.

Consider the stratification of $\BB_{\F_\ell}$ by Schubert cells.  For
a dominant weight $\mu$ Kempf showed that the global Cousin complex
$C\bul_{\F_\ell}(\mu)$, $$ C^{k}_{\F_\ell}(\mu)=\underset{w\in
W,\lth(w)=k}{\bigoplus}H^k_{S_{w,\F_\ell}}(\BB_{\F_\ell},\L(\mu)), $$
provides a resolution of the contragradient Weyl module $$\D
W_{\F_\ell}(\mu):=H^0(\BB_{\F_\ell},\L(\mu)).$$ Note also that the
algebra $U_{\F_\ell}(\g)$ maps naturally to the algebra of global
twisted differential operators on $\BB_{\F_\ell}$ with the twist
$\L(\mu)$ and thus it acts naturally on the above global Cousin
complex.  We claim that the obtained complex is isomorphic to the
specialization of the quasi-BGG complex $\D B_{\F_\ell}(\mu):=
\D B\bul_\CA(\mu)\ten_\CA\F_\ell$. Denote the quotient ring of 
$\CA$ by the ideal generated by the $\ell$-th cyclotomic polynomial 
by $\CA_\ell'$

It follows that $\D B_\ell(\mu)=\D B_{\CA_\ell'}(\mu)\ten \CC$ provides
a resolution of $\D W_\ell(\mu)$ since it can be viewed as a
specialization of the complex $\D B_{\CA_\ell'}(\mu)$ to the generic
point of $\Spec\CA_\ell'$, on the other hand by the above
considerations the specialization of the latter complex into $\Spec\
\F_\ell$ is exact.

Using the last result we prove the generalized Feigin conjecture.
Namely we show that $$ \Exts_{\u_\ell}(\underC,\D
M_\ell^w(w\cdot\ell\lambda))=
H_{T_{S_w}^*(\BB)}^{\sharp(R^+)}(\til{\CN},\pi^*\L(\lambda))$$ as a
module over both $U(\g)$ and
$H^0(\til{\CN},\O_{\til{\CN}})=H^0(\CN,\O_{\CN})$.  Here
$T_{S_w}^*(\BB)$ denotes the conormal bundle to the Schubert cell
$S_w$.  This statement allows us to calculate semiinfinite cohomology
of $\u_\ell$ with coefficients in the contragradient quasi-BGG
complex:  $$ \Exts_{\u_\ell}(\underC,\D B_\ell\bul(\ell\lambda))=
H_{\mu^{-1}(\n^+)}^{\sharp(R^+)}(\til{\CN},\pi^*\L(\lambda)).  $$

In the fourth section we formulate some further conjectures concerning
possible origin of the quasi-BGG complex.  A similar complex
expressing a Weyl module in terms of Verma modules exists in the
category $\O$ for the corresponding affine Lie algebra at a negative
rational level.  We hope that some extension of the Kazhdan-Lusztig
functor takes the latter complex to our quasi-BGG complex over
$\U_\ell$.

\subsection{Acknowledgements.}  The material of the paper is based
upon work supported by the U.S.  Civilian Research and Development
Foundation under Award No.  RM1-265 and by the grant INTAS 94-94720.

Part of the work over results of the paper was done during the
author's visit to Freiburg University in January 1998.  The author is
happy to thank the inviting orgainization for hospitality and
stimulating research conditions.

The author would like to thank A.  Braverman, R.  Bezrukavnikov, B.
Feigin, M.  Finkelberg, V.  Ginzburg, V.  Ostrik and W.  Soergel for
helpful discussions.

\section{Small quantum groups.}  We start with recalling notations and
the main results from \cite{Ar7}.

\subsection{Quantum groups at roots of unity.}  Fix a {\em Cartan
datum} $(I,\cdot)$ of the finite type and a {\em simply connected root
datum} $(Y,X,\ldots)$ of the type $(I,\cdot)$.  Thus we have
$Y=\Z[I]$, $X=\hom(Y,\Z)$, and the pairing $\langle\ ,\ \rangle:\
Y\times X\map\Z$ coincides with the natural one (see \cite{L3}, I 1.1,
I 2.2).  In particular the data contain canonical embeddings
$I\hookrightarrow Y, i\mapsto i$ and $I\hookrightarrow X, i\mapsto i':
\langle i',j\rangle:= 2i\cdot j/i\cdot i$.  The latter map is
naturally extended to an embedding $Y\subset X$.  Denote by $\hgt$ the
linear function on $X$ defined on elements $i', i\in I$, by
$\hgt(i')=1$ and extended to the whole $X$ by linearity.  The root
system (resp.  the positive root system) corresponding to the data
$(Y,X,\ldots)$ is denoted by $R$ (resp.  by $R^+$), below $W$ denotes
the Weyl group of $R$.

Like in \cite{Ar7}, Section 2, we denote the {\em Drinfeld-Jimbo
quantum group} defined over the field $\Q(v)$ of rational functions in
$v$ (resp.  the {\em Lusztig version of the quantum group} defined
over $\CA=\Z[v,v^{-1}]$) by $\U$ (resp.  by $\U_{\CA}$).  Fix an odd
number $\ell$ satisfying the conditions from \cite{GK} and a primitive
$\ell$-th root of unity $\xi$.  Define a $\CC$-algebra
$\til{\U}_\ell:=\U_{\CA}\ten_{\CA}\CC$, where $v$ acts on $\CC$ by
multiplication by $\xi$.  It is known that the elements $K_i^{\ell},
i\in I$, are central in $\til{\U}_\ell$.  Set
$\U_\ell:=\til{\U}_\ell/(K_i^\ell-1, i\in I)$.  The algebra $\U_\ell$
is generated by the elements $E_i$, $E_i^{(\ell)}$, $F_i$,
$F_i^{(\ell)}$, $K_i^{\pm1}$, $i\in I$.  Here $E_i^{(\ell)}$ (resp.
$F_i^{(\ell)}$) denotes the {\em $\ell$-th quantum divided power} of
the element $E_i$ (resp.  $F_i$) specialized at the root of unity
$\xi$.

\subsubsection{Quantum groups in positive charavteristic.} Like in 
\cite{L1} consider the specialization of $\U_\CA$ in characteristic $\ell$. 
Namely let $\CA_\ell'$ be the quotient of $\CA$ by the ideal generated by the
 $\ell$-th cyclotomic polynomial. T
hen $\CA_\ell'/(v-1)$ is isomorphic to the finite field ${\mathbb F}_\ell$.
 Thus the algebraic closture $\F_\ell$
becomes a $\CA$-algebra. We set $\U_{\F_\ell}:=\U_\CA\ten_{\CA}\F_\ell$. 
Ir is known that the elements $K_i$, $i\in I$, are central in 
$\U_{\F_\ell}$ and the algebra $\U_{\F_\ell}/(K_i,i\in I)$ is isomorphic to 
$U_{\Z}(\g)\ten \F_\ell$, where $U_\Z(\g)$ denotes the Kostant integral form 
for the unicersal enveloping algebra of 
$\g$. 

\sssn
Following Lusztig we define the {\em small quantum group} $\u_\ell$ at
the $\ell$-th root of unity $\xi$ as the subalgebra in $\U_\ell$ generated by
all $E_i, F_i, K_i^{\pm1}, i\in I$.  Denote the subalgebra in
$\u_\ell$ generated by $E_i, i\in I$ (resp.  $F_i, i\in I$, resp.
$K_i, i\in I$), by $\u_\ell^+$ (resp.  $\u_\ell^-$, resp.
$\u_\ell^0$).  The subalgebra $\u_\ell^-\ten\u_\ell^0$ (resp.
$\u_\ell^0\ten\u_\ell^+$) in $\u_\ell$ is denoted by $\b_\ell^-$
(resp.  by $\b_\ell^+$).

Recall that a finite dimensional algebra $A$ is called {\em Frobenius}
if the left $A$-modules $A$ and $A^*:=\hom_{\CC}(A,\CC)$ are
isomorphic.

\sssn \Lemma (see \cite{Ar1}, Lemma 2.4.5) \label{one} The algebras
$\u_\ell^+$ and $\u_\ell^-$ are Frobenius.\qed

Note that the algebra $\u_\ell$ is graded naturally by the abelian
group $X$.  Using the function $\hgt$ we obtain a $\Z$-grading on
$\u_\ell$ from this $X$-grading.  In particular the subalgebra
$\u_\ell^+$ (resp.  $\u_\ell^-$) is graded by $\Z_{\ge0}$ (resp.  by
$\Z_{\le0}$).

Below we present several well known facts about the algebra $\u_\ell$
to be used later.  Recall that an augmented subalgebra $B\subset A$
with the augmentation ideal $\overline{B}\subset B$ is called {\em
normal} if $A\overline{B}=\overline{B}A$.  If so, the space
$A/A\overline{B}$ becomes an algebra.  It is denoted by $A//B$.  Fix
an augmentation on $\u_\ell$ as follows:  $E_i\mapsto0, F_i\mapsto 0,
K_i\mapsto 1$ for every $i\in I$.  Set
$\underC:=\u_\ell/\overline{\u}_\ell$.

\sssn \Lemma (see \cite{AJS} 1.3, \cite{L2} Theorem 8.10)
\label{mainu} \begin{itemize} \item[(i)] The multiplication in
$\u_\ell$ provides a vector space isomorphism
$\u_\ell=\u_\ell^-\ten\u_\ell^0\ten\u_\ell^+$; the subalgebra
$\u^0_\ell$ is isomorphic to the group algebra of the group
$(\Z/\ell\Z)^{\sharp(I)}$.  \item[(ii)] The subalgebra $\u_\ell\subset
\U_\ell$ is normal and we have $\U_\ell//\u_\ell=U(\g)$.\qed
\end{itemize}

Denote the category of $X$-graded finite dimensional left
$\u_\ell$-modules $M=\underset{\lambda\in X}{\bigoplus}M_\lambda$ such
that $K_i$ acts on $M_\lambda$ by multiplication by the scalar
$\xi^{\langle i,\lambda\rangle}$ and $E_i:\ M_\lambda\map
M_{\lambda+i'},\ F_i:\ M_\lambda\map M_{\lambda-i'}$ for all $i\in I$,
with morphisms preserving $X$-gradings, by $\u_\ell\mod$.  For $M,N\in
\u_\ell\mod$ and $\lambda\in\ell\cdot X$ we define the shifted module
$M\langle\lambda\rangle\in\u_\ell\mod:\
M\langle\lambda\rangle_\mu:=M_{\lambda+\mu}$ and set
$\hom_{\u_\ell}(M,N):=\underset{\lambda\in\ell\cdot
X}{\bigoplus}\hom_{\u_\ell\mod}(M\langle \lambda\rangle,N)$.
Evidently the spaces $\hom_{\u_\ell}(\cdot,\cdot)$ and posess natural
$\ell\cdot X$-gradings.

\subsection{Cohomology of small quantum groups.}  Consider the
$\ell\cdot X\times\Z$-graded algebra
$\Ext_{\u_\ell}(\underC,\underC)$.  Note that by Shapiro lemma and
Lemma~\ref{mainu} (ii) the Lie algebra $\g$ acts naturally on the
$\ext$ algebra and the multiplication in the algebra satisfies Lebnitz
rule with respect to the $\g$-action.  In \cite{GK} Ginzburg and Kumar
obtained a nice description of the multiplication structure as well as
the $\g$-module structure on $\Ext_{\u_\ell}(\underC,\underC)$ as
follows.

\subsubsection{Functions on the nilpotent cone.}  Let $G$ be the
simply connected Lie group with the Lie algebra $\g$.  Then $G$ acts
on $\g$ by adjunction.  The action preserves the set of nilpotent
elements $\CN\subset\g$ called the {\em nilpotent cone} of $\g$.  The
action is algebraic, thus it provides a morphism of $\g$ into the Lie
algebra of algebraic vector fields on the nilpotent cone ${\mathcal
V}ect(\CN)$.  The latter algebra acts on the algebraic functions
$H^0(\CN,\O_{\CN})$.  The action is $G$-integrable.  Note also that
the natural action of the group $\CC^*$ provides a grading on
$H^0(\CN,\O_{\CN})$ preserved by the $G$-action.

\sssn \Theorem (see \cite{GK}) The algebra and $\g$-module structures
on $H^0(\CN,\O_{\CN})$ and on $\Ext_{\u_\ell}(\underC,\underC)$
coincide.  The homological grading on the latter algebra corresponds
to the grading on the former one provided by the $\CC^*$-action.  The
$X$-grading on $H^0(\CN,\O_{\CN})$ provided by the weight
decomposition with respect to the action of the Cartan subalgebra in
$\g$ corresponds to the natural $\ell\cdot X$-grading on the space
$\Exts_{\u_\ell}(\underC,\underC)$.  \qed

\subsection{Semiinfinite cohomology of small quantum
groups.}\label{setup} Here we present the definition of semiinfinite
cohomology of a finite dimensional associative algebra $\u$.  The
setup for the definition includes a nonnegatively (resp.
nonpositively) graded subalgebras $\u^+$ and $\u^-$ in $\u$ such that
the multiplication in $\u$ provides a vector space isomorphism
$\u^-\ten\u^+\til{\map}\u$.  Below we suppose that the algebra $\u^+$
is Frobenius.  Note still that the general definition of associative
algebra semiinfinite cohomology given in \cite{Ar1} requires neither
this restriction nor the assumption that $\dim \u<\infty$.

Consider first the {\em semiregular} $\u$-bimodule
$S_\u^{\u^+}=\u\ten_{\u^+}\u^{+*}$, with the right $\u$-module
structure provided by the isomorphism of left $\u^+$-modules
$\u^+\cong\u^{+*}$.  Note that $S_{\u}^{\u^+}$ is free over the
algebra $\u$ both as a right and as a left module.

For a complex of graded $\u$-modules $M\bul=\underset{p,q\in\Z}
{\bigoplus}M_p^q,\ d:\ M_p^q\map M_p^{q+1}$ we define the support of
$M\bul$ by $\supp M\bul:=\{(p,q)\in\Z^2|M_p^q\ne0\}$.  We say that a
complex $M\bul$ is {\em concave} (resp.  {\em convex}) if there exist
$s_1,s_2\in\N, t_1,t_2\in\Z$ such that $\supp
M\bul\subset\{(p,q)\in\Z^2| s_1q+p\le t_1, s_2q-p\le t_2\}$ (resp.
$\supp M\bul\subset\{(p,q)\in\Z^2| s_1q+p\ge t_1, s_2q-p\ge t_2\}$).

Let $M\bul,N\bul\in{\mathcal C}om(\u\mod)$.  Suppose that $M\bul$ is
convex and $N\bul$ is concave.  Choose a convex (resp.  concave)
complex $R_{\uparrow}\bul(M\bul)$ (resp.  $R_{\downarrow}\bul(N\bul)$)
in ${\mathcal C}om(\u\mod)$ quasiisomorphic to $M\bul$ (resp.
$N\bul$) and consisting of $\u^+$-projective (resp.
$\u^-$-projective) modules.

\sssn \Def We set $$ \Exts_{\u}(M\bul,N\bul)
:=H\bul(\Hom_{\u}(R\bul_{\uparrow}(M\bul),S_{\u}^{\u^+}\ten_\u
R_{\downarrow}\bul(N\bul))).  $$ \sssn \Lemma (see.  [Ar1] Lemma
3.4.2, Theorem 5) The spaces $\Exts_{\u}(M\bul,N\bul)$ do not depend
on the choice of resolutions and define functors $$
\ext_{\u}^{\si+k}(\cdot,\cdot):\ \u\mod\times\u\mod\map{\mathcal
V}ect,\ k\in\Z.\qed $$ Below we consider semiinfinite cohomology of
algebras $\u_\ell$, $\b^+_\ell$, $\b_\ell^-$ etc.  with coefficients
in $X$-graded modules.  The $\Z$-grading on such a module is obtained
from the $X$-grading using the function $\hgt:\ X\map \Z$.

Evidently the spaces $\Exts_{\u_{\ell}}(M\bul,N\bul)$ posess natural
 $\ell\cdot X$-gradings.  The following statement is a direct
 consequence of Lemma~\ref{mainu}(ii).

\sssn \Lemma \label{action} Let $M,N\in\u_\ell\mod$ be restrictions of
some $\U_\ell$-modules.  Then the spaces $\Exts_{\u_{\ell}}(M,N)$ have
natural structures of $\g$-modules, and the $\ell\cdot X$-gradings on
them coincide with the $X$-gradings provided by the weight
decompositions of the modules with respect to the standard Cartan
subalgebra in $\g$.\qed

\subsubsection{Semiinfinite cohomology of the trivial
${\u}_\ell$-module.}  Note that for a $\U_\ell$-module $M$ the space
$\Exts_{\u_{\ell}}(\underC,M)$ has a natural structure of a module
over the algebra $\Ext_{\u_{\ell}}(\underC,\underC)$ and this
structure is equivariant with respect to the action of the Lie algebra
$\g$.

Let us recall the geometric description of semiinfinite cohomology of
the trivial $\u_\ell$-module obtained in \cite{Ar7}.  Consider the
standard positive nilpotent subalgebra $\n^+\subset\g$ as a Zarisski
closed subset in $\CN$.  Consider the $H^0(\CN,\O_{\CN})$-module
$H_{\n^+}^{\sharp(R^+)}(\CN,\O_{\CN})$ of local cohomology with
supports in $\n^+\subset\CN$ for the structure sheaf $\O_\CN$.  Recall
also that the space $H_{\n^+}^{\sharp(R^+)}(\CN,\O_{\CN})$ has a
natural $\g$-module structure defined as follows.  First the Lie
algebra ${\mathcal{V}}ect(\CN)$ acts naturally on the local cohomology
space.  Now the adjoint action of $G$ on $\CN$ defines a Lie algebra
inclusion $\g\subset {\mathcal{V}}ect(\CN)$.  Note that the
$H^0(\CN,\O_{\CN})$-module $H_{\n^+}^{\sharp(R^+)}(\CN,\O_{\CN})$ is
evidently $\g$-equivariant with respect to the defined actions of $\g$
on the algebra and the module.

Note also that the subset $\n^+\subset \CN$ is $\CC^*$-stable, thus
the space $H_{\n^+}^{\sharp(R^+)}(\CN,\O_{\CN})$ is naturally graded
by the $\CC^*$-action.

The following statement sums up the main results from \cite{Ar7}.
\vskip 1mm \noindent \Theorem \label{main} \begin{itemize} \item[(i)]
Semiinfinite cohomology of the trivial $\u_\ell$-module vanishes in
even degrees.  \item[(ii)] The space
$\ext_{\u_\ell}^{\si+2k+1}(\underC,\underC)$ is isomorphic to the
homogeneous component in $H_{\n^+}^{\sharp(R^+)}(\CN,\O_{\CN})$ of the
weight $k$.  \item[(iii)] $H_{\n^+}^{\sharp(R^+)}(\CN,\O_{\CN})$ is
isomorphic to $\Exts_{\u_{\ell}}(\underC,\underC)$ both as a module
over $\Ext_{\u_{\ell}}(\underC,\underC)=H^0(\CN,\O_{\CN})$ and as a
$\g$-module.\qed \end{itemize}

\subsubsection{Springer-Grothendieck resolution of the nilpotent
cone.}  We provide another description of the space
$\Exts_{\u_{\ell}}(\underC,\underC)$ in terms of local cohomology as
follows.

Choose a maximal torus $H$ in the simply connected Lie group $G$ with
the Lie algebra $\g$.  The choice provides the root decomposition of
$\g$ and in particular its triangular decomposition
$\g=\n^-\oplus\h\oplus\n^+$.  Consider the Borel subgroup $B\subset G$
with the Lie algebra $\b^+=\h\oplus\n^+$ and the flag variety $\BB$.
The group $G$ acts on $\BB$ by left translations and the restriction
of this action to $B$ is known to have finitely many orbits.  These
orbits are isomorphic to affine spaces and called {\em the Schubert
cells}.  The Bruhat decomposition of $G$ shows that the Schubert cells
are enumerated by the Weyl group.  Denote the orbit corresponding to
the element $w\in W$ by $S_w$.

It is well known that the cotangent bundle $T^*(\BB)=:\til{\CN}$ has a
nice realization $\til{\CN}=\{(B_x,n)|n\in \operatorname{Lie}(B_x)\}
$, where $B_x$ denotes some Borel subgroup in $G$ and $n$ is a
nilpotent element in the Lie algebra $\operatorname{Lie}(B_x)$.  The
map $$ \mu:\ \til{\CN}\map \CN,\ (B_x,n)\mapsto n, $$ is known to be a
resolution of singularities of $\CN$ called {\em the
Springer-Grothendieck resolution}.  Note that the map $\mu$ is
$G$-equivariant.

Recall the following statement from \cite{Ar1}.

\sssn \Prop (see e.  g.  \cite{CG}, 3.1.36) \begin{itemize} \item[(i)]
$\mu^{-1}(\n^+)=\underset{w\in W}{\bigsqcup}T^*_{S_w}(\BB)$, where
$T_{S_w}^*(\BB)$ denotes the conormal bundle to $S_w$ in $\BB$.
\item[(ii)] $H^0(\til{\CN},\O_{\til{\CN}})=H^0(\CN,\O_{\CN})$ and the
higher cohomology spaces of the structure sheaf on $\til{\CN}$ vanish.
\item[(iii)] $H^{\sharp(R^+)}_{\n^+}(\CN,\O_{\CN})\til{\map}
H^{\sharp(R^+)}_{\mu^{-1}(\n^+)}(\til{\CN},\O_{\til{\CN}})$ both as a
$\g$-module and as a $H^0(\CN,\O_{\CN})$-module.  \qed \end{itemize}
\Cor $\Exts_{\u_{\ell}}(\underC,\underC)$ is isomorphic to
$H^{\sharp(R^+)}_{\mu^{-1}(\n^+)}(\til{\CN},\O_{\til{\CN}})$ both as a
module over
$\Ext_{\u_{\ell}}(\underC,\underC)=H^0(\til{\CN},\O_{\til{\CN}})$ and
as a $\g$-module.\qed

\subsection{Semiinfinite cohomology of contragradient Weyl modules.}
Our main goal now is to find a nice geometric interpretation of
semiinfinite cohomology of some remarkable $\u_\ell$-modules that
would generalize the results of stated in the previous section.

Fix the natural triangular decompositions of the algebra $\U_{\CA}$
(resp.  $\U_\ell$) as follows:
$\U_{\CA}=\U_{\CA}^-\ten\U_{\CA}^0\ten\U_{\CA}^+$ (resp.
$\U_{\ell}=\U_{\ell}^-\ten\U_{\ell}^0\ten\U_{\ell}^+$), where the
positive (resp.  negative) subalgebras are generated by the quantum
divided powers of the positive (resp.  negative) root generators in
the corresponding algebra.  We call the subalgebra
$\U_{\CA}^+\ten\U_{\CA}^0$ (resp.  $\U_{\ell}^+\ten\U_{\ell}^0$) the
{\em positive quantum Borel subalgebra} in $\U_{\CA}$ (resp.  in
$\U_\ell$) and denote it by $\B_{\CA}^+$ (resp.  by $\B^+_\ell$).  The
negative Borel subalgebras $\B_{\CA}^-$ and $\B_\ell^-$ are defined in
a similar way.

\ssn Fix a dominant integral weight $\lambda\in X$.  Consider the
module over the ``big'' quantum group $\U_{\CA}$ given by $$ \D
W_{\CA}(\lambda):=\left(\Coind_{\B^-_{\CA}}^{\U_{\CA}}\CC(\lambda)
\right)^{\operatorname{fin}} \text{ (resp.  by }
W_{\CA}(\lambda):=\left(\Ind_{\B^+_\ell}^{\U_\ell}\CC(\lambda)
\right)_{\operatorname{fin }}).  $$ Here $(*)^{\operatorname{fin}}$
(resp.  $(*)_{\operatorname{fin}}$) denotes the maximal finite
dimensional submodule (resp.  quotient module) in $(*)$.  The module
$\D W_{\CA}(\lambda)$ (resp.  $W_{\CA}(\lambda)$) is called {\em the
contragradient Weyl module} (resp.  {\em the Weyl module}) over
$\U_{\CA}$ with the highest weight $\lambda$.  We denote by $\D
W_\ell(\lambda)$ (resp.  by $W_\ell(\lambda)$) the specializations of
the corresponding modules into the chosen $\ell$-th root of unity.

It is known that at the generic value $\xi$ of the quantizing
parameter $v$ both modules $W_{\CA}(\lambda)$ and $\D
W_{\CA}(\lambda)$ specialize into the finite dimensional simple module
$L(\lambda)$ over the quantum group $\U_\xi:=\U_{\CA}\ten_{\CA}\CC$.
In particular we have $$\ch(W_{\CA}(\lambda))=\ch(\D
W_{\CA}(\lambda))=\sum_{w\in
W}{\frac{e^{w\cdot\lambda}}{\prod_{\alpha\in R^+}(1-e^\alpha)}}, $$
just like in the Weyl character formula in the semisimple Lie algebra
case.  Note also that $W_{\CA}(0)=\D W_{\CA}(0)=\underC$.

Below we consider semiinfinite cohomology of the algebra $\u_\ell$
with coefficients in the contragradient Weyl module with a
$\ell$-divisible highest weight $\ell\lambda$.  Ou considerations were
motivated by the following result of Ginzburg and Kumar (see
\cite{GK}).

Let $p$ denote the projection $\til{\CN}\map \BB$.  Consider the
linear bundle $\L(\lambda)$ on $\BB$ with the first Chern class equal
to $\lambda\in X=H^2(\BB,\Z)$.

\sssn \Theorem \label{Weyl} \begin{itemize} \item[(i)]
$\ext_{\u_\ell}^{\operatorname{odd}}(\underC,\D
W_{\ell}(\ell\lambda))=0$; \item[(ii)]
$\ext_{\u_\ell}^{2\bullet}(\underC,\D W_\ell(\ell\lambda))=
H^0(\til{\CN},p^*\L(\lambda))$ as a $H^0(\CN,\O_{\CN})$-module.  \qed
\end{itemize}

The following conjecture provides a natural semiinfinite analogue for
Theorem~\ref{Weyl}.

\sssn \label{maincon} \Con $\ext_{\u_\ell}^{\si+\bullet}(\underC,\D
W_{\ell}(\ell\lambda))=
H^0_{\mu^{-1}(\n^+)}(\til{\CN},p^*\L(\lambda))$ as a
$H^0(\CN,\O_{\CN})$-module.  The homological grading on the left hand
side of the equality corresponds to the grading by the natural action
of the group $\CC^*$ on the right hand side.  \qed

\sssn \Cor \label{chformula} \begin{gather*} \ch\left(
\ext_{\u_\ell}^{\si+\bullet}(\underC,\D
W_{\ell}(\ell\lambda)),t\right)\\=
\frac{t^{-\sharp(R^+)}}{\prod_{\alpha\in
R^+}(1-e^{-\ell\alpha})}\sum_{w\in W}
\frac{e^{w(\ell\lambda)}t^{2l(w)}}{\prod_{\alpha\in R^+,w(\alpha)\in
R^+}(1-t^2e^{-\ell\alpha}) \prod_{\alpha\in R^+,w(\alpha)\in
R^-}(1-t^{-2}e^{-\ell\alpha})}.\qed \end{gather*} Below we present the
main steps for the proof of the conjecture.  Details of the proof will
be given in the fourthcoming paper \cite{Ar8}.

\subsection{Local cohomology with coefficients in $p^*\L(\lambda)$.}
For a quasicoherent sheaf $\M$ on $\til{\CN}$ consider the natural
pairing $ H^i_{\mu^{-1}(\n^+)}(\til{\CN},\O_{\til{\CN}})\times
H^0(\til{\CN},\M) \map H^i_{\mu^{-1}(\n^+)}(\til{\CN},\M).  $
Evidently it is equivariant with respect to the
$H^0(\til{\CN},\O_{\til{\CN}})$-action.  Thus we obtain a
$H^0(\til{\CN},\O_{\til{\CN}})$-module morphism $$ s:\
H_{\mu^{-1}(\n^+)}^i(\til{\CN},\O_{\til{\CN}})\ten_{H^0(\til{\CN},
\O_{\til{\CN}})}
H^0(\til{\CN},\M)\map H^i_{\mu^{-1}(\n^+)}(\til{\CN},\M).  $$

\sssn \Prop \label{step1} For $\M=p^*\L(\lambda)$ the map $s$ is an
isomorphism.  \qed

\subsubsection{Similar construction for semiinfinite cohomology.}  We
will need some more homological algebra.  Fix a graded algebra $A$
with a subalgebra $B\subset A$.  Recall that in \cite{V} and
\cite{Ar1} the notion of a complex of graded $A$-modules K-semijective
with respect to the subalgebra $B$ was developed.  The following
statement gives an analogue of the standard technique of projective
resolutions in the semiinfinite case.

\sssn \Theorem (see \cite{Ar1}, Appendix B) Let
$SS_{\u_\ell^+}\bul(*)$ (resp.  $SS_{\u_\ell^-}\bul(*)$) denote a
K-semijective concave (resp.  convex) resolution of the
$\u_\ell$-module $(*)$ with respect to the subalgebra $\u_\ell^+$
(resp.  $\u_\ell^-$).  Then for a finite dimensional graded
$\u_\ell$-module $M$ we have \begin{itemize} \item[(i)]
$H\bul(\Hom_{\u_\ell}( SS\bul_{\u_\ell^-}(\underC),
SS\bul_{\u_\ell^-}(M))=\Ext_{\u_\ell}(\underC,M)$; \item[(ii)]
$H\bul(\Hom_{\u_\ell}( SS\bul_{\u_\ell^+}(\underC),
SS\bul_{\u_\ell^-}(\underC))=\Exts_{\u_\ell}(\underC,M)$.  \qed
\end{itemize} \sssn \Cor The composition of morphisms provides a
natural pairing $$
\ext_{\u_\ell}^{\si+i}(\underC,\underC)\times\ext^j_{\u_\ell}
(\underC,M)\map \ext^{\si+i+j}_{\u_\ell}(\underC,M).\qed $$ In
particular we obtain a $\Ext_{\u_\ell}(\underC,\underC)$-module map $$
\Exts_{\u_\ell}(\underC,\underC)\ten_{\Ext_{\u_\ell}(\underC,\underC)}
\Ext_{\u_\ell}(\underC,M)\map \Exts_{\u_\ell}(\underC,M).  $$
Combining this construction with the previous considerations we obtain
the following statement.

\sssn \Prop There exists a natural $H^0(\CN,\O_{\CN})$-module morphism
$$ \sigma:\
H_{\mu^{-1}(\n^+)}^{\sharp(R^+)}(\til{\CN},p^*\L(\lambda))\map
\Exts_{\u_\ell}(\underC,\D W(\ell\lambda)).  $$ \dok Follows from
Proposition~\ref{step1}.  \qed

Below we show that the morphism $\sigma$ is an isomorphism.  The main
tool for the demonstration of this fact is the {\em quasi-BGG complex}
providing a specialization of the classical BGG resolution for a
finite dimensional simple $\U$-module $L(\ell\lambda)$ into the root
of unity $\xi$.  This complex constructed below consists of direct
sums of $\U_\ell$-modules called {\em the quasi-Verma modules}.  On
the other hand we show that semiinfinite cohomology with coefficients
in quasi-Verma modules has a nice geometrical interpretation.

\section{The construction of the quasi-BGG complex.}  Recall that the
usual Bernshtein-Gelfand-Gelfand resolution of a finite dimensional
representatiln $L(\lambda)$ of the simple Lie algebra $\g$ has a nice
geometric interpretation as follows.  First by Borel-Weyl-Bott theorem
we have $L(\lambda)=H^)(\CB,\L(\lambda))$.  Next consider the
sratification $\{S_w,w\in W\}$ of the flag variety by the Schubert
cells.  Kempf showed that the contragradient BGG-resolution of
$L(\lambda)$ coincides with the global Cousin complex
$\C\bul(\lambda)$:  $$ \C^k(\lambda)=\underset{w\in
W,\lth(w)=k}{\bigoplus}H^k_{S_w}(\CB,\L(\lambda)).  $$ In particular
this construction inedtifies the local cohomology space
$H^k_{S_w}(\CB,\L(\lambda))$ with the contragradient Verma module $\D
M(w\cdot\lambda)$.

\subsection{Cousin complexes in positive characteristic.} \label{cousin} 
 In fact the
local cohomology construction due to Kempf works in a more general
setting.  Let $\F_\ell$ be the algebraic closture of the finite field
of characteristic $\ell$.  It is known that both the flag variety for
$\g$ and its stratification by Schubert cells are well defined over
$\F_\ell$.  We denote the corresponding objectsby $\CB_{\F_\ell}$ and
$S_{w,\F_\ell}$ respectively.

Let $U_{\F_\ell}(\g)$ be the specialization of the Kostant integral
form for the universal enveloping algebra $U_{\Z}(\g)$ into $\F_\ell$.
Fix a dominant weight $\lambda$.  Then it is known that just like in
the complex case $U_{\F_\ell}(\g)$ maps naturally into the algebra of
global differential operators on $\CB_{\F_\ell}$ with coefficients in
the line bundle $\L(\lambda)$ denoted by
$\operatorname{Diff}(\L(\lambda))$.

Consider the contragradient Weyl module over $U_{\F_\ell}(\g)$ defined
as follows:  $$ \D
W_{\F_\ell}(\lambda):=\left(\Coind_{U_{\F_\ell}(\b^-)}^{U_{\F_\ell}(\g)}
\F_\ell(\lambda)\right)^{\operatorname{fin}}.
$$ The Borel-Weyl-Bott theorem remains true in prime characteristic.

\vskip 1mm \noindent \Theorem \begin{itemize} \item[(i)]
$H^{>0}(\CB_{\F_\ell},\L(\lambda))=0$; \item[(ii)]
$H^{0}(\CB_{\F_\ell},\L(\lambda))=\D W(\lambda)$.\qed \end{itemize}

What is even more important for us is that the Kempf construction of
the global Cousin complex works over $\F_\ell$ as well.

\sssn \Theorem (see \cite{K}) For any dominant weight $\lambda$ there
exists a complex $\C_{\F_\ell}\bul(\lambda)$, $$
\C_{\F_\ell}^k(\lambda):=\underset{w\in
W,\lth(w)=k}{\bigoplus}H^k_{S_{w,\F_\ell}}(\CB_{\F_\ell},\L(\lambda)).
$$ \Rem The algebra $\operatorname{Diff}(\L(\lambda))$ acts naturally
on $\C_{\F_\ell}\bul(\lambda)$ and this complex becomes the one of
$U_{\F_\ell}(\g)$-modules.  An important difference from the complex
case is that the complex no longer consists of direct sums of
contragradient Verma modules.

Our aim is to mimick algebraically the above construction of the
global Cousin complex in the setting of the quantum group $\U_\ell$
rather than the algebra $U_{\F_\ell}(\g)$.

\subsection{Twisted quantum parabolic subalgebras in $\U_\CA$} Recall
that Lusztig has constructed an action of the {\em braid group}
$\mathfrak B$ corresponding to the Cartan data $(I,\cdot)$ by
automorphisms of the quantum group $\U_{\CA}$ well defined with
respect to the $X$-gradings (see \cite{L1}, Theorem 3.2).  Fix a
reduced expression of the maximal length element $w_0\in W$ via the
simple reflection elements:  $$ w_0=s_{i_1}\ldots
s_{i_{\sharp(R^+)}},\ i_k\in I.  $$ Then it is known that this reduced
expression provides reduced expressions for all the elements $w\in W$:
$ w=s_{i_1^w}\ldots s_{i_{l(w)}^w},\ i_k\in I.  $

Consider the standard generators $\{T_i\}_{i\in I}$ in the braid group
$\mathfrak B$.  Lifting the reduced expressions for the elements $w$
from $W$ into $\mathfrak B$ we obtain the set of elements in the braid
group of the form $T_w:= T_{i_1^w}\ldots T_{i_{l(w)}^w}$.

In particular we obtain the set of {\em twisted Borel subalgebras}
$w(\B_\CA^+)=T_w(\B_\CA^+)\subset\U_\CA$.  Note that
$w_0(\B_\CA^+)=\B_\CA^-=\U_\CA^-\ten\U_\CA^0$.

Fix a subset $J\subset I$ and consider the {\em quantum parabolic
subalgebra} $\P_{J,\CA}\subset\U_\CA$.  By definition this subalgebra
in $\U_\CA$ is generated over $\U_\CA^0$ by the elements $E_i$, $i\in
I$, $F_j$, $j\in J$, and by their quantum divided powers.  The
previous construction provides the set of {\em twisted quantum
parabolic subalgebras} $w(\P_{J,\CA}):=T_w(\P_{J,\CA})$ of the type
$J$ with the twists $w\in W$.

Note that the triangular decomposition of the algebra $\U_\CA$
provides the ones for the algebras $w(\B_\CA^+)$ and $w(\P_{J,\CA})$:
$$ w(\B_\CA^+) = (w(\B_\CA^+))^- \ten \U_\CA^0 \ten (w(\B_\CA^+))^+
\text{ and } w(\P_{J,\CA}) = (w(\P_{J,\CA}))^- \ten \U_\CA^0 \ten
(w(\P_{J,\CA})^+, $$ where $\left(w(\B_\CA^+)\right)^+ =
w(\B_\CA^+)\cap \U_\CA^-$, $\left(w(\P_{J,\CA})\right)^+ =
w(\P_{J,\CA})\cap \U_\CA^-$ etc.

Choosing the root of unity $\xi$ and specializing the algebra $\CA$
into $\CC$, where $v$ acts on $\CC$ by the multiplication by $\xi$, we
obtain in particular the subalgebras $w(\B_\ell^+)\subset\U_\ell$,
$w(\P_{J,\ell})\subset \U_\ell$, $w(\b_\ell^+)\subset\u_\ell$,
$w(\p_{J,\ell})\subset \u_\ell$ with the induced triangular
decompositions.

\subsection{Semiinfinite induction and coinduction} From now on we
will use freely the technique of associative algebra semiinfinite
homology and cohomology for a graded associative algebra $A$ with two
subalgebras $B,N\subset A$ equipped with a triangular decomposition
$A=B\ten N$ on the level of graded vector spaces.  We will not recall
the construction of these functors referring the reader to \cite{Ar1}
and \cite{Ar2}.

Let us mention only that these functors are bifunctors ${\mathsf
D}(A\mod)\times{\mathsf D}(\oppA\mod)\map{\mathsf D}({\mathcal V}ect)$
where the associative algebra $\oppA$ is defined as follows.

Consider the semiregular $A$-module $S_A^N:=A\ten_NN^*$.  It is proved
in \cite{Ar2} that under very weak conditions on the algebra $A$ the
module $S_A^N$ is isomorphic to the $A$-module
$\left(S_A^N\right)':=\hom_B(A,B)$.  Thus $\End_A(S_A^N)\supset
N^\opp$ and $\End_A(S_A^N)\supset B^\opp$ as subalgebras.  The algebra
$\oppA$ is defined as the subalgebra in $\End_A(S_A^N)$ generated by
$B^\opp$ and $N^\opp$.  It is proved in \cite{Ar2} that the algebra
$\oppA$ has a triangular decomposition $\oppA=N^\opp\ten B^\opp$ on
the level of graded vector spaces.  Yet for an arbitrary algebra $A$
the algebras $\oppA$ and $A^\opp$ do not coincide.

However the following statement shows that in the case of quantum
groups that correspond to the root data $(Y,X,\ldots)$ of the {\em
finite} type $(I,\cdot)$.  the equality of $A^\opp$ and $\oppA$ holds.

\sssn \Prop We have \begin{itemize} \item[(i)] $\U^{\sharp}=\U^\opp$,
$\U_\CA^{\sharp}=\U_\CA^\opp$, $\U_\ell^{\sharp}=\U_\ell^\opp$;
\item[(ii)] $w(\B^+)^\sharp= w(\B^+)^\opp$, $w(\B_\CA^+)^\sharp=
w(\B_\CA^+)^\opp$, $w(\B_\ell^+)^\sharp= w(\B_\ell^+)^\opp$,
$w(\P_{J,\CA})^\sharp= w(\P_{J,\CA})^\opp$, $w(\P_{J,\ell})^\sharp=
w(\P_{J,\ell})^\opp$.  \end{itemize} \dok The first part is proved
similarly to Lemma 9.4.1 from \cite{Ar6}.  The second one follows
immediately from the first one.  \qed

\sssn \Def Let $M\bul$ be a convex complex of $w(\B_\CA^+)$-modules.
By definition set \begin{gather*}
\SInd_{w(\B_\CA^+)}^{\U_\CA}(M\bul):=
\tor_{\si+0}^{w(\B_\CA^+)}(S_{\U_\CA}^{\U_\CA^+},M\bul) \text{ and }\\
\SCoind_{w(\B_\CA^+)}^{\U_\CA}(M\bul):=
\ext^{\si+0}_{w(\B_\CA^+)}(S_{\U_\CA}^{\U_\CA^-},M\bul).
\end{gather*} The functors $\SInd_{w(\B_\ell^+)}^{\U_\ell}(\cdot)$,
$\SCoind_{w(\B_\ell^+)}^{\U_\ell}(\cdot)$,
$\SInd_{w(\P_{J,\ell}^+)}^{\U_\ell}(\cdot)$,
$\SCoind_{w(\P_{J,\ell}^+)}^{\U_\ell}(\cdot)$ etc.  are defined in a
similar way.

\sssn \Lemma (see \cite{Ar4}) \begin{itemize} \item[(i)]
$\tor_{\si+k}^{w(\B_\CA^+)}(S_{\U_\CA}^{\U_\CA^+},\cdot)=0$ for
$k\ne0$; \item[(ii)]
$\ext^{\si+k}_{w(\B_\CA^+)}(S_{\U_\CA}^{\U_\CA^-},\cdot)=0$ for
$k\ne0$; \item[(iii)] $\SInd_{w(\B_\CA^+)}^{\U_\CA}(\cdot)$ and
$\SCoind_{w(\B_\CA^+)}^{\U_\CA}(\cdot)$ define {\em exact} functors
$w(\B_\CA^+)\mod\map\U_\CA\mod$.\qed \end{itemize} Similar statements
hold for the algebras $w(\B_\ell^+)$, $w(\P_{J,\CA})$ and
$w(\P_{J,\ell})$.

\subsection{Quasi-Verma modules} We define the {\em quasi-Verma
module} over the algebra $\U_\CA$ (resp.  $\U_\ell$) with the highest
weight $w\cdot\lambda$ by $$ M_\CA^w(w\cdot \lambda):=
\SInd_{w(\B_\CA^+)}^{\U_\CA}(\CA(\lambda)) \text{ (resp.  }
M_\ell^w(w\cdot \lambda):=
\SInd_{w(\B_\ell^+)}^{\U_\ell}(\CC(\lambda))).  $$ The {\em
contragradient quasi-Verma module} $ \D M_\CA^w(w\cdot\lambda)$ (resp.
$ \D M_\ell^w(w\cdot\lambda)$) is defined by $$ \D M_\CA^w(w\cdot\lambda):=
\SCoind_{w(\B_\CA^+)}^{\U_\CA}(\CA(\lambda)) \text{ (resp.  } \D
M_\ell^w(w\cdot \lambda):= \SCoind_{w(\B_\ell^+)}^{\U_\ell}(\CC(\lambda))).
$$ We list the main properties of quasi-Verma modules.

\sssn \Prop (see \cite{Ar4}) \begin{itemize} \item[(i)] Fix a dominant
integral weight $\lambda\in X$.  Suppose that $\xi\in\CC^*$ is not a
root of unity.  Then the $\U_\xi$-module
$M_\xi^w(w\cdot\lambda):=M_\CA^w(w\cdot\lambda)\ten_{\CA}\CC$ (resp.
$\D M_\xi^w(w\cdot\lambda):=\D M_\CA^w(w\cdot\lambda)\ten_{\CA}\CC$)
is isomorphic to the usual Verma module $M_\xi(w\cdot\lambda)$ (resp.
to the usual contragradient Verma module $\D M_\xi(w\cdot\lambda)$).
\item[(ii)] For any $\lambda\in X$ we have \begin{gather*}
\ch(M_\CA^w(w\cdot\lambda)) = \ch(\D M_\CA^w(w\cdot\lambda)) =
\ch(M_\ell^w(w\cdot\lambda)) \\ = \ch(\D M_\ell^w(w\cdot\lambda)) =
\frac{e^{w\cdot\lambda}}{\prod_{\alpha\in R^+}(1-e^{-\alpha})}.\qed
\end{gather*} \end{itemize} Thus for a dominant weight $\lambda$ one
can consider $M_\CA^w(w\cdot\lambda)$ as a flat family of modules over
the quantum group for various values of the quantizing parameter with
the fiber at the generic point equal to the Verma module
$M_\xi(w\cdot\lambda)$.  \vskip 1mm \noindent \Rem Note that by
definition $M_\ell^e(\lambda)=M_\ell(\lambda)$ and $\D
M_\ell^e(\lambda)=\D M_\ell(\lambda)$, where $e$ denotes the unity
element of the Weyl group.  In particular for a dominant weight
$\lambda$ we have a natural projection $M_\ell^e(\lambda)\map
W_\ell(\lambda)$ and a natural inclusion $\D
W_\ell(\lambda)\hookrightarrow \D M_\ell^e(\lambda)$.

\subsection{The $\U_\ell({\frak sl}_2)$ case.}  Let us investigate throughly
quasi-Verma modules in the case of $\U_\ell=\U_\ell({\frak sl}_2)$.
First we find the simple subquotient modules in the module
$M_\ell^e(k\ell)=M_\ell(k\ell).$

Recall the classification of the simple objects in the category of
$X$-graded $\U_\ell^0$-semisimple $\U_\ell$-modules locally finite
with respect to the action of $E_i$ and $E_i^{(\ell)}$, $i\in I$,
obtained by Lusztig in \cite{L2}.  In the ${\frak sl}_2$ case it looks
as follows.  Identify the weight lattice $X$ with $\Z$.  \vskip 1mm
\noindent \Prop \begin{itemize} \item[(i)] For $0\le k<\ell$ the
simple $\u_\ell({\frak sl}_2)$-module $L(k)$ is a restriction of a
simple $\U_\ell({\frak sl}_2)$-module.  \item[(ii)] Any simple
$\U_\ell({\frak sl}_2)$-module from the category described above is
isomorphic to a module of the form $L(k)\ten L(m\ell)$, where $0\le
k<\ell$.  Here the simple module $L(m\ell)$ is obtained from the
simple $U({\frak sl}_2)$-module $L(m)$ via restricton using the map
$\U_\ell({\frak sl}_2)\map\U({\frak sl}_2)//\u_\ell({\frak
sl}_2)=U({\frak sl}_2)$.  \qed \end{itemize} Denote the only
reflection in the Weyl group for ${\frak sl}_2$ by $s$.  \vskip 1mm
\noindent \Lemma \begin{itemize} \item[(i)] $\ch M_\ell^e(0)=\ch
L(0)+\ch L(-2)+\ch L(-2\ell)$.  \item[(ii)] For $k>0$ we have $\ch
M_\ell^e(k\ell)=\ch L(k\ell)+\ch L(k\ell-2)+\ch L(-k\ell-2)+\ch
L(-(k+2)\ell)$.  \item[(iii)] $\ch M_\ell^s(s\cdot 0)=\ch L(-2)+\ch
L(-2\ell)$.  \item[(iv)] For $k>0$ we have $\ch M_\ell^s(s\cdot
k\ell)=+\ch L(-k\ell-2)+\ch L(-(k+2)\ell)$.  \qed \end{itemize} In
fact it is easy to find the filtrations on quasi-Verma modules with
simple subquotients that correspond to the character equalities above.
\vskip 1mm \noindent \Lemma \begin{itemize} \item[(i)] For $k>0$ there
exist exact sequences \begin{gather*} 0\map L(k\ell-2)\map
W(k\ell)\map L(k\ell)\map 0,\\ 0\map L(-k\ell-2)\map M_\ell^s(s\cdot
k\ell)\map L(-(k+2)\ell)\map 0.  \end{gather*} \item[(ii)] There
exists a filtration on $M_\ell^e(0)$ with subquotients as follows:  $$
\gr^1M_\ell^e(0)=L(0),\ \gr^2M_\ell^e(0)=L(-2),\
\gr^3M_\ell^e(0)=L(-2\ell).  $$ \item[(iii)] For $k>0$ there exists a
filtration on $M_\ell^e(k\ell)$ with subquotients as follows:
\begin{gather*} \gr^1M_\ell^e(k\ell)=L(k\ell),\
\gr^2M_\ell^e(k\ell)=L(k\ell-2),\\
\gr^3M_\ell^e(k\ell)=L(-(k+2)\ell),\
\gr^4M_\ell^e(k\ell)=L(-k\ell-2).\qed \end{gather*} \end{itemize} Thus
we obtain the following statement.

\sssn \label{sl2case} \Prop For any $k\ge0$ there exists an exact
sequence of $\U_\ell=\U_\ell({\frak sl}_2)$-modules $$ 0\map
M_\ell^s(s\cdot k\ell)\map M_\ell^e(k\ell)\map W_\ell(k\ell)\map0.\qed $$

A more accurate calculation of characters shows that the complex similar 
to the one from the previous Proposition exists also for non-$\ell$-divisible
 dominant highest weights.
\vskip 1mm
\noindent
 \Prop For any $\mu\ge0$ there exists an exact
sequence of $\U_\ell=\U_\ell({\frak sl}_2)$-modules $$ 0\map
M_\ell^s(s\cdot \mu)\map M_\ell^e(\mu)\map W_\ell(\mu)\map0.\qed $$

\Rem \begin{itemize} \item[(i)] In fact it is easy to verify that for $\mu$ 
dominant the module
$M_\ell^s(s\cdot \mu)$ is isomorphic to the {\em contragradient}
Verma module $\D M_\ell(s\cdot \mu)$.  \item[(ii)] Note that if
$\xi$ is not a root of unity then the usual BGG resolution in the
${\frak sl}_2$ case provides an exact complex $$ 0\map M_\xi^s(s\cdot
k\ell)\map M_\xi^e(k\ell)\map L(k\ell)\map0.  $$ \end{itemize} Thus we
see that the flat family of such complexes over
$\CC^*\setminus\{$roots of unity$\}$ is extended over the whole
$\CC^*$.
\subsubsection{The $\U_\CA({\frak sl}_2)$-case.}
In fact we will need the existence of {\em a complex of 
$\U_\CA({\frak sl}_2)$-modules} 
similar to the one constructed in the previous subsection. Recall that by 
definition there exists a natural projection of $\U_{\CA}$-modules 
$M_{\CA}(\mu)\map W_{\CA}(\mu)$. Denote its kernel by $K$. 

\sssn
\Prop
The $\U_{\CA}({\mathfrak {sl}}_2)$-module $K$ is isomorphic to 
$\D M_{\CA}(s\cdot\mu)$.

\dok
First note that $K$ has a highest weight vector $p$ of the weight 
$s\cdot\mu$. In particular there exists a natural map 
$K\map\D M_{\CA}(s\cdot\mu)$.
Thus it is sufficient to check that the module $K$ is co-generated by this 
vector. This claim follows from the calculation below. Denote the highest 
weight vector in $M_{\CA}(\mu)$ by $p_\mu$. Then evidently 
$p=F^{(\mu+1)}p_{\mu}$. Using the formulas from
~\cite{L1}, 6.4, for every $m>\mu$ we have
$$
E^{(m-\mu-1)}\cdot F^{(m)}p_\mu=F^{(\mu+1)}\prod_{s=1}^{m-\mu-1}\frac{Kv^{-\mu
-1-s+1}-K^{-1}v^{\mu+1+s-1}}{v^s-s^{-s}}v_\mu=-F^{(\mu+1)}p_\mu.
\qed
$$
\Cor 
For every positive integer $\mu$ there exists an exact complex of $\U_{\CA}
({\mathfrak {sl}}_2)$-modules
$$ 0\map
M_\CA^s(s\cdot \mu)\map M_\CA^e(\mu)\map W_\CA(\mu)\map0.\qed $$
We call the complex $M_\CA^s(s\cdot \mu)\map M_\CA^e(\mu)$ the quasi-BGG 
complex for the weight $\mu$ and denote it by $B_\CA\bul(\mu)$.

\subsubsection{The $U_{\F_\ell}({\frak sl}_2)$-case.} The crucial point in 
the proof of the exactness of the quasi-BGG complex for  arbitrary root data 
$(Y,X,\ldots)$  uses a geometric argument in positive characteristic. Thus 
the following Lemma is impor
tant being the $\frak {sl}_2$ case of the general picture.
\vskip 1mm
\noindent
\Lemma \label{sl2ref}
\begin{itemize}
\item[(i)] For each $i\in I$ the action of $K_i-1$ on the complex of
 $\U_{\F_\ell}$-modules $B_{\F_\ell}\bul(\mu):=B_\CA\bul(\mu)\ten_\CA\F_\ell$
 is trivial, i.~e. $B_{\F_\ell}\bul(\mu)$ becomes a complex of $U_{\F_\ell}
(\g)$-modules.
\item[(ii)] The complex $\D B_{\F_\ell}\bul(\mu)$ is isomorphic to the global 
Cousin complex
$C\bul(\mu)$
for the line bundle $\L(\mu)$  on $\P_{\F_\ell}^1$  discussed 
in~\ref{cousin}.\qed
\end{itemize}
\subsection{Construction of the quasi-BGG complex.}  Here we extend
the previous considerations to the case of the quantum group $\U_\CA$
for arbitrary root data $(Y,X,\ldots)$ of the finite type $(I,\cdot)$.
Fix a dominant weight $\mu \in X$.

First we construct an inclusion
$M_\CA^{w'}(w'\cdot\mu)\hookrightarrow
M_\CA^{w}(w\cdot\mu)$ for a pair of elements $w',w\in W$ such
that $\lth(w')=\lth(w)+1$ and $w'>w$ in the Bruhat order on the Weyl group.
In fact we can do it explicitly only for $w'$ and $w$ differing by a
simple reflection:  $w'=ws_i$, $i\in I$.

Consider the twisted quantum parabolic subalgebra $w(\P_{i,\CA})$.
Then $w(\P_{i,\CA})\supset w(\B_\CA^+)$ and $w(\P_{i,\CA})\supset
ws_i(\B_\CA^+)$.  Consider also the Levi quotient algebra
$w(\P_{i,\CA})\map w(\LL_{i,\CA})$.  The algebra $w(\LL_{i,\CA})$
is isomorphic to $\U_\CA({\frak sl}_2)\ten_{\U_\CA^0({\frak
sl}_2)}\U_\CA^0$.

By Proposition \ref{sl2case} we have a natural inclusion of $w
(\LL_{i,\CA})$-modules
$$\SInd_{ws_i(\LL_{i,\CA}^+)}^{w(\LL_{i,\CA})}\CA(\mu)
\hookrightarrow
\SInd_{w(\LL_{i,\CA}^+)}^{w(\LL_{i,\CA})}\CA(\mu).$$ \Lemma
\begin{itemize} \item[(i)] $\SInd_{w(\B_\CA^+)}^{\U_\CA}(\cdot)=
\SInd_{w(\P_{i,\CA})}^{\U_\CA}\circ
\SInd_{w(\B_\CA^+)}^{w(\P_{i,\CA})}(\cdot)$.  \item[(ii)]
$\SInd_{ws_i(\B_\CA^+)}^{\U_\CA}(\cdot)=
\SInd_{w(\P_{i,\CA})}^{\U_\CA}\circ
\SInd_{ws_i(\B_\CA^+)}^{w(\P_{i,\CA})}(\cdot)$.  \item[(iii)]
$\SInd_{w(\B_\CA^+)}^{\U_\CA}(\CA(\mu))=
\SInd_{w(\P_{i,\CA})}^{\U_\CA}\circ
\Res_{w(\P_{i,\CA})}^{w(\LL_{i,\CA})}\circ
\SInd_{W(\LL_{i,\CA}^+)}^{w(\LL_{i,\CA})}(\mu)$.
\item[(iv)] $\SInd_{ws_i(\B_\CA^+)}^{\U_\CA}(\CA(\mu))=
\SInd_{w(\P_{i,\CA})}^{\U_\CA}\circ
\Res_{w(\P_{i,\CA})}^{w(\LL_{i,\CA})}\circ
\SInd_{ws_i(\LL_{i,\CA}^+)}^{w(\LL_{i,\CA})}(\mu)$.  \qed
\end{itemize} \noindent \Cor For $w'=ws_i>w$ in the Bruhat order we
have a natural inclusion of $\U_\CA$-modules $i_\CA^{ws_i,w}:\
M_\CA^{ws_i}(ws_i\cdot\mu)\hookrightarrow
M_\CA^w(w\cdot\mu).  \qed $

Recall that if $v$ acts on $\CC$ by 
$\xi$ that is  not a root of unity then the $\U_\xi$-module 
$M_\xi^w(w\cdot\mu):=M_\CA^w(w\cdot\mu)\ten_\CA\CC$ is
isomorphic to the usual Verma module $M_\xi(w\cdot\mu)$.  Thus
the morphism $i_\xi^{w',w}$ coincides with the standard inclusion of
Verma modules constructed by J.  Bernshtein, I.M.  Gelfand and S.I.
Gelfand in \cite{BGG} that becomes a component of the differential in
the BGG resolution.  In other words we see that the flat family of
inclusions $i_\xi^{ws_i,w}:\
M_\ell^{ws_i}(ws_i\cdot\mu)\hookrightarrow
M_\xi^w(w\cdot\mu) $ defined   for  $\xi\in\CC^*\setminus\{$roots of
unity$\}$ can be extended naturally over the whole $\Spec \CA$.

Iterating the inclusion maps we obtain a flat family of submodules
$i^w_\xi(M_\xi^w(w\cdot\mu))\subset M_\xi^e(\mu)$ for
$\xi\in\Spec \CA$, $w\in W$, providing an extension of the standard
lattice of Verma submodules in $M_\xi(\mu)$ defined a priori for
$\xi\in\CC^*\setminus\{$roots of unity$\}$.

\sssn \Lemma For a pair of elements $w',w\in W$ such that
$\lth(w')=\lth (w)+1$ and $w'>w$ in the Bruhat order we have $$
i_\CA^{w'}(M_\CA^{w'}(w'\cdot\mu)) \hookrightarrow
i_\CA^w(M_\CA^{w}(w\cdot\mu)).  $$ \dok To prove the
statement note that the condition $\{A_\xi$ is a submodule in
$B_\xi\}$ is a closed condition in a flat family.  \qed

Now using the standard combinatorics of the classical BGG resolution
we obtain the following statement.

\sssn \Theorem There exists a complex of $\U_\CA$-modules
$B\bul_CA(\mu)$ with $$
B_\CA^{-k}(\mu)=\underset{w\in W,\lth(w)=k}{\bigoplus}
M_\CA^w(w\cdot\mu) $$ and with differentials provided by
direct sums of the inclusions $i_\CA^{w',w}$.  \qed

\sssn \Def We call the complex $B\bul_\CA(\mu)$ the {\em
quasi-BGG complex} for the dominant weight $\mu\in X$.

Denote the complex of $\u_\ell$-modules  $B_\CA(\mu)\ten_\CA\CC$
 by $B_\ell\bul(\mu)$. Below we prove that for $\ell$ prime 
the quasi-BGG complex  $B_\ell\bul(\mu)$ is in fact quasiisomorphic to the 
Weyl module 
$W_\ell(\mu)$.

\subsection{Quasi-BGG complexes in positive characteristic.} Note that we 
can perform specialization of the quasi-BGG complex over $\CA$ into a root 
of unity in two steps. Consider the complex of $\U_{\CA'_\ell}$-modules
$\D B_{\CA_\ell'}\bul(\mu) :=\D B_{\CA}\bul(\mu)\ten_\CA\CA_\ell'$. 
Evidently its specialization into the generic point of $\Spec \CA_\ell'$ 
coincides with $\D B_{\ell}\bul(\mu)$.  On the other hand consider the 
specialization of the complex into the spec
ial point $\Spec\F_\ell\hookrightarrow\Spec\CA_\ell'$. The following 
statement is similar to Lemma~\ref{sl2ref}.

\sssn
\Prop
\begin{itemize}
\item[(i)]
 For each $i\in I$ the action of $K_i-1$ on the complex of 
$\U_{\F_\ell}$-modules $B_{\F_\ell}\bul(\mu):=B_{\CA_\ell'}\bul(\mu)
\ten_{\CA_\ell'}\F_\ell$ is trivial, i.~e. $B_{\F_\ell}\bul(\mu)$ becomes a 
complex of $U_{\F_\ell}(\g)$-modules.
\item[(ii)] The complex  is isomorphic to the global Cousin complex
$C\bul(\mu)$
for the line bundle $\L(\mu)$  on $\CB_{\F_\ell}$  discussed in~\ref{cousin}.
\qed
\end{itemize}
\Cor 
The complex $\D B_{\F_\ell}\bul(\mu)$ is  quasiisomorphic to the 
contragradient Weyl module $\D W_{\F_\ell}(\mu)=\D W_{\CA}(\mu)\ten_\CA\F_\ell
=H^0(\CB_{\F_\ell},\L(\mu))$.\qed
\vskip 1mm
\noindent
Surprisingly this result proves the exactness of the quasi-BGG complex over 
$\U_\ell$.
\vskip 1mm
\noindent
\Theorem \label{qis}
For $\ell$ prime the  complex $\D B_{\F_\ell}\bul(\mu)$ is  quasiisomorphic 
to the contragradient Weyl module $\D W_{\ell}(\mu)$.

\dok 
Denote the evident morphism $\D W_{\CA_\ell'}(\mu)\map\D B_{\CA_\ell'}\bul
(\mu)$ by ${can}_{\CA_\ell'}$ . Consider the complex  $\Cone\left(
{can}_{\CA_\ell'}\right)$. By the previous Corollary the specialization of  
the complex into the special point $\Spec\F_\ell\hookrightarrow\Spec 
\CA_\ell'$ is exact. By the Nakayama lemma it follows that the specialization 
of the complex into the generic point of $\Spec \CA_\ell'$ is also exact. In 
particular the complex $\Cone\left({can}_{\ell}\right):=\Cone\left({ca
n}_{\CA_\ell'}\right)\ten_{\CA_\ell'}\CC$, where $v$ acts on $\CC$ by the 
$\ell$-th root of unity, is exact as well. \qed

\subsection{Semiinfinite cohomology with coefficients in quasi-Verma
modules.}  Recall the following construction that plays crucial role
in considerations of Ginzburg and Kumar in \cite{GK}.

Let $\left(\B_\ell^+\mod\right)^{\operatorname{fin}}$ (resp.
$\left(\U_\ell\mod\right)^{\operatorname{fin}}$, resp.
$\left(U(\b^+)\mod\right)^{\operatorname{fin}}$, resp.
$\left(U(\g)\mod\right)^{\operatorname{fin}}$) be the category of
finite dimensional $X$-graded modules over the corresponding algebra
with the action of the Cartan subalgebra semisimple and well defined
with respect to the $X$-gradings.  Consider the functors:
\begin{gather*} \Coind_{\B_\ell^+}^{\U_\ell}:\
\B_\ell^+\mod\map\U_\ell\mod;\
\left(\Coind_{\B_\ell^+}^{\U_\ell}\right)^{\operatorname{fin}}:\
\left(\B_\ell^+\mod\right)^{\operatorname{fin}}\map\left(\U_\ell
\mod\right)^{ \operatorname{fin}};\\ \Coind_{U(\b^+)}^{U(\g)}:\
U(\b^+)\mod\map U(\g)\mod;\\
\left(\Coind_{U(\b^+)}^{U(\g)}\right)^{\operatorname{fin}}:\
\left(U(\b^+)\mod\right)^{\operatorname{fin}}\map\left(U(\g)\mod
\right)^{ \operatorname{fin}},\\ (\cdot)^{\b_\ell^+}:\
\B_\ell^+\mod\map U(\b^+)\mod \text{ and }
\left(\B_\ell^+\mod\right)^{\operatorname{fin}}\map
\left(U(\b^+)\mod\right)^{\operatorname{fin}};\\ (\cdot)^{\u_\ell}:\
\U_\ell\mod\map U(\g)\mod \text{ and }
\left(\U_\ell\mod\right)^{\operatorname{fin}}\map
\left(U(\g)\mod\right)^{\operatorname{fin}}, \end{gather*} where
$(\cdot)^{\operatorname{fin}}$ denotes taking the maximal finite
dimensional submodule in $(\cdot)$ and $(\cdot)^{\b_\ell^+}$ (resp.
$(\cdot)^{\u_\ell}$) denotes taking $\b_\ell^+$- (resp.
$\u_\ell$)-invariants.  \vskip 1mm \noindent \Prop (see \cite{GK})
\begin{itemize} \item[(i)] $(\cdot)^{\u_\ell}\circ
\Coind_{\B_\ell^+}^{\U_\ell}= \Coind_{U(\b^+)}^{U(\g)}=
(\cdot)^{\b_\ell^+}$; \item[(ii)] $(\cdot)^{\u_\ell}\circ
\left(\Coind_{\B_\ell^+}^{\U_\ell}\right)^{\operatorname{fin}}=
\left(\Coind_{U(\b^+)}^{U(\g)}\right)^{\operatorname{fin}}
(\cdot)^{\b_\ell^+}$.\qed \end{itemize} The semiinfinite analogue for
the first part of the previous statement looks as follows.  Fix $w\in
W$.  Consider the functors:  \begin{gather*}
\SCoind_{w(\B_\ell^+)}^{\U_\ell}:\ {\mathsf
D}(w(\B_\ell^+)\mod)\map{\mathsf D} (\U_\ell\mod),\\
\SCoind_{U(w(\b^+))}^{U(\g)}:\ {\mathsf D}(U(w(\b^+))\mod)\map{\mathsf
D} (U(\g)\mod),\\ \Exts_{w(\b_\ell^+)}(\underC,\cdot):\ {\mathsf
D}(w(\B_\ell^+)\mod)\map{\mathsf D} (U(w(\b^+))\mod),\\
\Exts_{\u_\ell}(\underC,\cdot):\ {\mathsf
D}(w(\U_\ell)\mod)\map{\mathsf D} (U(\g))\mod).  \end{gather*} \sssn
\Theorem We have $$ \Exts_{\u_\ell}(\underC,\cdot) \circ
\SCoind_{w(\B_\ell^+)}^{\U_\ell} = \SCoind_{U(w(\b^+))}^{U(\g)} \circ
\Exts_{w(\b_\ell^+)}(\underC,\cdot).  $$ \dok To simplify the
notations we work with semiinfinite homology and semiinfinite
induction instead of semiinfinite cohomology and semiinfinite
coinduction.  By \cite{Ar1}, Appendix B, every convex complex of
$w(\B_\ell^+)$-modules is quasiisomorphic to a K-semijective convex
complex.  Consider a K-semijective convex complex of
$w(\B_\ell^+)$-modules $SS\bul$.  Note that both semiinfinite
induction functors are exact and take K-semijective complexes to
K-semijective complexes.  Thus we have \begin{gather*}
\Tors^{\u_\ell}(\underC,\cdot)\circ
\SInd_{w(\B_\ell^+)}^{\U_\ell}(SS\bul)= \Tors^{\u_\ell}(\underC,
\Tors^{w(\B_\ell^+)}(S_{\U_\ell}^{\U_\ell^+},SS\bul))\\=
\Tors^{w(\B_\ell^+)}( \Tors^{\u_\ell}(\underC,
S_{\U_\ell}^{\U_\ell^+}),SS\bul))= \Tors^{w(\B_\ell^+)}
(S_{U(\g)}^{U(\n^+)},SS\bul)\\= \Tors^{w(\B_\ell^+)} (\underC,
S_{U(\g)}^{U(\n^+)})\ten SS\bul)= \Tors^{U(w(\b^+))}( \underC,
\Tors^{w(\b_\ell^+)}(\underC, S_{U(\g)}^{U(\n^+)}) \ten SS\bul))\\=
\Tors^{U(w(\b^+))}( S_{U(\g)}^{U(\n^+)}), \Tors^{w(\b_\ell^+)}(
\underC,SS\bul)) = \SInd_{U(w(\b^+))}^{U(\g)} \circ
\Tors^{w(\b_\ell^+)}(\underC,SS\bul).  \end{gather*} Here we used the
fact that the subalgebra $w(\b_\ell^+)\subset w(\B_\ell^+)$ is normal
with the quotient algebra equal to $U(w(\b^+))$.  \qed \vskip 1mm
\noindent \Cor $$ \Exts_{\u_\ell}(\underC,\D
M_\ell^w(w\cdot\ell\lambda))=
H_{T_{S_w}^*(\BB)}^{\sharp(R^+)}(\til{\CN},\pi^*\L(\lambda))$$ as a
module over both $U(\g)$ and
$H^0(\til{\CN},\O_{\til{\CN}})=H^0(\CN,\O_{\CN})$.  \qed

Recall that for the Springer-Grothendieck resolution of the nilpotent
cone $\mu:\ \til{\CN}\map \CN$ we have $\mu^{-1}(\n^+)=\underset{w\in
W}{\bigsqcup}T^*_{S_w}(\BB)$.  \vskip 1mm \noindent \Prop
\begin{itemize} \item[(i)] There exists a filtration on
$H^{\sharp(R^+)}_{\mu^{-1}(\n^+)}(\til{\CN},\pi^*\L(\lambda))$ with
the subquotients equal to
$H^{\sharp(R^+)}_{T^*_{S_w}(\BB)}(\til{\CN},\pi^*\L(\lambda))$, for
$w\in W$.  \item[(ii)] $ \Exts_{\u_\ell}(\underC,\D
B_\ell\bul(\ell\lambda))=
H_{\mu^{-1}(\n^+)}^{\sharp(R^+)}(\til{\CN},\pi^*\L(\lambda))$.
 \qed \end{itemize}
Comparing this statement with Theorem~\ref{qis} we obtain the main
result of the section.
\sssn \Theorem For $\ell$ prime
$\ext_{\u_\ell}^{\si+\bullet}(\underC,\D W(\ell\lambda))=
H^0_{\mu^{-1}(\n^+)}(\til{\CN},p^*\L(\lambda))$ both as a $\U(\g)$-module and 
as a $H^0(\til{\CN},\O_{\til{\CN}})$-module.
\qed

\section{Further results and conjectures.}  In this section we present
several facts without proof.  We also formulate some conjectures
concerning possible origin of the quasi-BGG complex.

\subsection{Alternative triangular decompositions of $\u_\ell$.}
\label{alt} Note that the definition of semiinfinite cohomology starts
with specifying a triangular decomposition of a {\em graded} algebra
$\u$.  Fix a subset $J\subset I$.  Instead of the usual height
function consider the linear map $\hgt_J:\ X\map\Z$ defined on the
elements $i', i\in I$ by $\hgt_J(i')=0$ for $i\in J$ and
$\hgt_J(i')=1$ otherwize and extended to the whole $X$ by linearity.
Now we work in the category of complexes of $X$-graded
$\u_\ell$-modules satisfying conditions of concavity and convexity
with respect to the $\Z$-grading obtained from the $X$-grading with
the help of the function $\hgt_J$.

Consider the triangular decomposition of the small quantum group
$\u_\ell=\p_{J,\ell}^-\ten \u_{J,\ell}^+$, where $\p_{J,\ell}^-$
denotes the small quantum negative parabolic subalgebra in $\u_\ell$
corresponding to the subset $J\subset I$ and $\u_{J,\ell}^+$ denotes
the quantum analogue of the nilpotent radical in $\p_J^+$ defined with
the help of Lusztig generators of $\U_\ell$ and $\u_\ell$ (see
\cite{L1}).

Then it is known that the subalgebra $\u_{J,\ell}^+$ in $\u_\ell$ is
Frobenius just like $\u_\ell^+$.  Thus it is possible to use the
general definition of semiinfinite cohomology presented in
\ref{setup}.  Denote the corresponding functor by
$\Exts_{\u_\ell,J}(*,*)$.

On the other hand consider the classical negative parabolic subalgebra
$\p_J^-\subset\g$ and its nillradical $\n_J^-$.  Choose the standard
$X$-homogeneous root basis $\{f_\alpha\}$ in the space $\n_J^-$.
Consider the subset in $\CN^{(J)}\subset\CN$ annihilated by all the
elements of the base dual to $\{f_\alpha\}$.

\sssn \Theorem $\Exts_{\u_\ell,J}(\underC,\underC)\til{\map}
H^{\dim(\n^-_J)}_{\CN^{(J)}}(\CN,\O_{\CN})$ as
$H^0(\CN,\O_{\CN})$-modules.  \qed

\subsection{Contragradient Weyl modules with non-$\ell$-divisible
highest weights.}  Fix a {\em dominant} weight of the form
$\ell\lambda+ w\cdot 0$, where $\lambda\in X$ and $w\in W$.  It is
known that all the dominant weights in the linkage class containing
$0$ look like this.  Consider the contragradient Weyl module $\D
W(\ell\lambda+ w\cdot0)$.  The following statement generalizes
Corollary \ref{chformula}.  Its proof is similar to the proof of
Conjecture \ref{maincon}.

\sssn \Theorem \begin{gather*} \ch\left(
\ext_{\u_\ell}^{\si+\bullet}(\underC,\D W(\ell\lambda+
w\cdot0)),t\right)\\= \frac{t^{-\sharp(R^+)+l(w)}}{\prod_{\alpha\in
R^+}(1-e^{-\ell\alpha})}\sum_{v\in W}
\frac{e^{v(\ell\lambda)}t^{2l(v)}}{\prod_{\alpha\in R^+,v(\alpha)\in
R^+}(1-t^2e^{-\ell\alpha}) \prod_{\alpha\in R^+,v(\alpha)\in
R^-}(1-t^{-2}e^{-\ell\alpha})}.\qed \end{gather*}
\subsection{Contragradient Weyl modules:  alternative triangular
decompositions.}  Fix the triangular decomposition of the small
quantum group $\u_\ell$ like in \ref{alt}.  A natural generalization
of Conjecture~\ref{maincon} to the case of the parabolic triangular
decomposition looks as follows.  We keep the notations from~\ref{alt}.

\sssn \Con $\Exts_{\u_\ell,J}(\underC,\D W(\ell\lambda))\til{\map}
H^{\dim(\n^-_J)}_{\CN^{(J)}}(\CN,\mu_*p^*\L(\lambda))$ as a module
over $H^0(\CN,\O_{\CN})$.  \qed

\subsection{Connection with affine Kac-Moody algebras} Finally we
would like to say a few words about a possible explanation for the
existence of quasi-Verma modules and quasi-BGG resolutions.

Suppose for simplicity that the root data $(Y,X,\ldots)$ are {\em
simply laced}, i.  e.  the corresponding Cartan matrix is symmetric.
Consider the affine Lie algebra $\hat\g=\g\ten\CC[t,t^{-1}]\oplus\CC
K$ corresponding to $\g$.  Fix a {\em negative} level $-h^\vee+k$,
where $k\in 1/2\Z_{<0}\setminus\Z_{<0}$ and $h^\vee$ denotes the dual
Coxeter number for chosen root data of the finite type.  Consider the
Kazhdan-Lusztig category $\til{\mathcal O}_{-k}$ of
$\g\ten\CC[t]$-integrable finitely generated $\hat\g$-modules
diagonizible with respect to the Cartan subalgebra in $\hat\g$ at the
level $-2h^\vee-k$.  Kazhdan and Lusztig showed that the category
$\til{\mathcal O}_k$ posasses a structure of a rigid tensor category
with the {\em fusion} tensor product $\overset{\cdot}{\ten}$.
Moreover, they proved the following statement.

\sssn \Theorem (see \cite{KL}) Let $\ell=-2k$.  Then the tensor
category $(\til{\mathcal O}_k,\overset{\cdot}{\ten})$ is equivalent to
the category $\left(\U_\ell\mod\right)^{\operatorname{fin}}$ with the
tensor product provided by the Hopf algebra structure on $\U_\ell$.
\qed

We denote the functor $ (\til{\mathcal O}_k,\overset{\cdot}{\ten})\map
\left(\left(\U_\ell\mod\right)^{\operatorname{fin}},\ten\right)$
providing the equivalence of categories by $\til{\mathsf {kl}}$.
Consider the usual category ${\mathcal O}_k$ for $\hat\g$ at the same
level.  Finkelberg constructed a functor ${\mathsf {kl}}:\ {\mathcal
O}_k\map\U_\ell\mod$ extending the functor $\til{\mathsf {kl}}$ (see
\cite{F}).  Note that the functor $\mathsf{kl}$ has no chance to be an
equivalence of categories because it is known not to be exact.

Fix a dominant (resp.  arbitrary) weight $\lambda\in X$.  Consider now
the contragradient Weyl module $\D{\mathcal
W}(\lambda)=\Coind_{U(\g\ten\CC[t])}^{U_{-2h^\vee-k}(\hat\g)}L(\lambda)$
and the contragradient Verma module $\D{\mathcal
M}(\lambda)=\Coind_{U(\g\ten\CC[t])}^{U_{-2h^\vee-k}(\hat\g)}\D
M(\lambda)$ over $\hat\g$ at the level $-2h^\vee-k$ , where
$L(\lambda)$ (resp.  $M(\lambda)$) denotes the simple module (resp.
the contragradient Verma module over $\g$ with the highest weight
$\lambda$.  Then the usual contragradient BGG resolution of
$L(\lambda)$ provides a resolution $\D\B(\lambda)$ of the
contragradient Weyl module $\D {\mathcal W}(\lambda)$ consisting of
direst sums of contragradient Verma modules of the form $\D{\mathcal
M}(w\cdot\lambda)$, where $\in W$.  It is known that the
Kazhdan-Lusztig functor takes Weyl and contragradient Weyl modules
over $\hat\g$ to Weyl (resp.  contragradient Weyl) modules over
$\U_\ell$.

\sssn \Con The functor ${\mathsf {kl}}$ takes $\D{\mathcal
M}(w\cdot\ell\lambda)$ to the contragradient quasi-Verma module $\D
M_\ell^w(w\cdot\ell\lambda)$.  Moreover the complex ${\mathsf
{kl}}(\D{\mathbf B}(\ell\lambda))$ is quasiisomorphic to $\D
W(\ell\lambda)$.\qed


\begin{thebibliography}{99} \bibitem[AJS]{AJS} H.H.Andersen,
J.C.Jantzen, W.Soergel.  {\it Representations of quantum groups at
p-th roots of unity and of semisimple groups in characteristic p:
independence of p.}  Asterisque \hbox{\bf 220}, (1994).
\bibitem[Ar1]{Ar1} S.M.Arkhipov.  {\it Semiinfinite cohomology of
quantum groups.}  Comm.  Math.  Phys.  \hbox{\bf 188} (1997), 379-405.
\bibitem[Ar2]{Ar2} S.M.Arkhipov.  {\it Semiinfinite cohomology of
associative algebras and bar duality.}  International Math.  Research
Notices No.  17 (1997), 833-863.  \bibitem[Ar3]{Ar3} S.M.Arkhipov.
{\it A new construction of the semi-infinite BGG resolution.}
Preprint q-alg/9605043 (1996), 1-29.  \bibitem[Ar4]{Ar4} S.M.Arkhipov.
{\it Semiinfinite cohomology of quantum groups II.}  Preprint q-alg
9610045 (1996), 1-42.  To appear in Amer.  Math.  Translations.
\bibitem[Ar5]{Ar5} S.M.Arkhipov.  {\it Semiinfinite cohomology of
small quantum groups.}  Funct.  An.  Appl., \hbox{\bf 32} No.1 (1998),
61-66.  \bibitem[Ar6]{Ar6} S.M.Arkhipov.  Ph.D.  Thesis, Moscow State
University (1997), 1-93.  \bibitem[Ar7]{Ar7} S.M.Arkhipov.  {\it A
proof of Feigin's conjecture.}  Mathematical Research Letters {\bf 5}
(1998), 1-20. 
\bibitem[Ar8]{Ar8}S.M.Arkhipov. {\it On a quantum analogue of the global 
Cousin complex for the flag manifold.} In preparation.
\bibitem[BGG]{BGG}I.~I.~Bernstein, I.~M.~Gelfand,
S.~I.~Gelfand.  {\it Differential operators on the principal affine space and
investigation of} $\g${\em -modules}, in Proceedings of Summer School on Lie
groups of Bolyai J\=anos Math. Soc., Helstead, New York, 1975.
\bibitem[CG]{CG} N.Chriss, V.Ginzburg.  {\it
Representation theory and complex geometry.}  Boston, Birkh\"auser,
(1995).  \bibitem[DCKP]{DCKP} C.De Concini, V.Kac, C.Procesi. {\it Some
remarkable degenerations of quantum groups.} Comm. Math. Phys.
\hbox{\bf 157}  (1993), p.405-427.

\bibitem[CG]{CG} N.Chriss, V.Ginzburg.  {\it
Representation theory and complex geometry.}  Boston, Birkh\"auser,
(1995).  
 \bibitem[F]{F} M.Finkelberg.  {\it An
equivalence of fusion categories.}  Harvard Ph.D.  thesis, (1993).
\bibitem[GeM]{GeM} S.I.Gelfand, Yu.I.Manin.  {\it Methods of
homological algebra.}  M:  Nauka, (1988).  (In Russian)
\bibitem[GK]{GK} V.Ginzburg, N.Kumar.  {\it Cohomology of quantum
groups at roots of unity.}  Duke Math.  J.  \hbox{\bf 69}, (1993),
179-198.  \bibitem[H]{H} W.H.Hesselink.  {\it On the character of the
nullcone.}  Math.  Ann.  \hbox{\bf 252}, (1980), 179-182.
\bibitem[K]{K} G.Kempf.  {\it The Grothendieck-Cousin complex of an
induced representation.}  Adv.  in Math., \hbox{\bf 29}, 310-396,
(1978).  \bibitem[KL 1,2,3,4]{KL} D.Kazhdan and G.Lusztig.  {\it
Tensor structures arising from affine Lie algebras.}  I, J.  Amer.
Math.  Soc.  {\bf 6}, (1993), 905-947; II, J.  Amer.  Math.  Soc.
{\bf 6} (1993), 949-1011; III, J.  Amer.  Math.  Soc.  {\bf 7} (1994),
335-381; IV, J.  Amer.  Math.  Soc.  {\bf 7}, (1994), 383-453.
\bibitem[L1]{L1} G.Lusztig.  {\it Quantum groups at roots of unity.}
Geom.  Dedicata {\bf 95}, (1990), 89-114.  \bibitem[L2]{L2} G.Lusztig.
{\it Modular representations and quantum groups.}  Contemp.  Math.
\hbox{\bf 82}, (1989), 59-77.  \bibitem[L3]{L3} G.Lusztig.  {\em
Introduction to quantum groups.}  (Progress in Mathematics \hbox{\bf
110}), Boston etc.  1993 (Birkh\"auser).  \bibitem[V]{V} A.Voronov.
{\it Semi-infinite homological algebra.}  Invent.  Math.  {\bf 113},
(1993), 103--146.  \end{thebibliography}
\end{document}